\documentclass[12pt]{amsart}
\usepackage{amssymb,latexsym,amsmath,amsthm,enumitem,mathrsfs,geometry}
\usepackage{fullpage}
\usepackage{fancyhdr}
\usepackage{xcolor}
\usepackage{bbm}
\usepackage{stmaryrd}
\usepackage{mathrsfs}
\usepackage{hyperref}
\usepackage{lineno}
\usepackage{diagbox}
\usepackage{array}
\usepackage{tikz}
\usetikzlibrary{arrows}
\usetikzlibrary{positioning}
\usepackage{float}
\usepackage{comment}
\usepackage{mathtools}
\usepackage{cleveref}

\newtheorem{theorem}{Theorem}[section]
\newtheorem*{theorem*}{Theorem}

\newtheorem{conjecture}{Conjecture}[section]
\newtheorem{lemma}{Lemma}[section]
\newtheorem*{remark*}{Remark}

\newtheorem{remark}{Remark}[section]

\newcommand{\R}{\mathbb{R}}
\newcommand{\N}{\mathbb{N}}

\begin{document}

\numberwithin{equation}{section}

\title{Rough numbers between consecutive primes}

\author[Gafni]{Ayla Gafni}
\address{Department of Mathematics, University of Mississippi, University MS 38677}
\email{argafni@olemiss.edu}
\author[Tao]{Terence Tao}
\address{Department of Mathematics, UCLA, 405 Hilgard Ave, Los Angeles CA 90024}
\email{tao@math.ucla.edu}

\subjclass[2020]{}
\thanks{}

\date{\today}

\begin{abstract} Using a sieve-theoretic argument, we show that almost all gaps $(p_n, p_{n+1})$ between consecutive primes $p_n, p_{n+1}$ contain a natural number $m$ whose least prime factor $p(m)$ is at least the length $p_{n+1} - p_n$ of the gap, confirming a prediction of Erd\H{o}s.  In fact the number $N(X)$ of exceptional gaps with $p_n \in [X,2X]$ is shown to be at most $O(X/\log^2 X)$.  Assuming a form of the Hardy--Littlewood prime tuples conjecture, we establish a more precise asymptotic $N(X) \sim c X / \log^2 X$ for an explicit constant $c>0$, which we believe to be between $2.7$ and $2.8$.  To obtain our results in their full strength we rely on the asymptotics for singular series developed by Montgomery and Soundararajan.
\end{abstract}

\maketitle

\section{Introduction and background} \label{intro}

\begin{figure}
\centering
\includegraphics[width=0.8\textwidth]{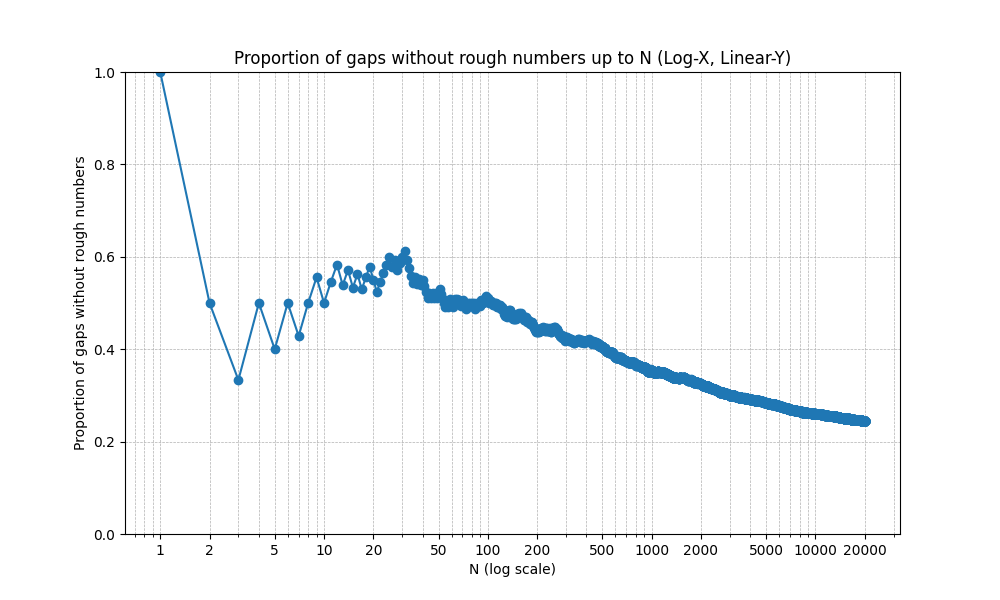}
\caption{Proportion of the first $N$ prime gaps that do not contain a rough number, for $N$ up to $2 \times 10^4$.  The horizontal axis is logarithmically scaled.  The code to generate this image (and the next) was initially generated by ChatGPT and then modified by the authors.}\label{gaps-fig}
\end{figure}

For a natural number $m$, let $p(m)$ denote the least prime factor of $m$ (with the convention $p(1)=1$).  We will informally refer to numbers $m$ with a large value of $p(m)$ as \emph{rough numbers}.  Let $p_n$ denote the $n^\mathrm{th}$ prime, thus the interval $(p_n, p_{n+1})$ is the $n^\mathrm{th}$ \emph{prime gap}; we refer to the length $p_{n+1} - p_n$ of this interval as the $n^\mathrm{th}$ \emph{prime gap length}.  Let us say that the $n^\mathrm{th}$ prime gap \emph{contains a rough number} if there exists $p_n < m < p_{n+1}$ such that
\begin{equation}\label{rough-gap}
    p(m) \geq p_{n+1} - p_n.
\end{equation}
For instance, any twin prime gap (for which $p_{n+1}-p_n=2$) will contain a rough number, since $p(p_n+1) \geq 2$.  The first few $n$ for which the $n^\mathrm{th}$ prime gap contains a rough number are
$$ 2, 3, 5, 7, 10, 13, 15, 17, 20, \dots.$$
Numerically, the proportion of $n$ for which the $n^\mathrm{th}$ prime gap does not contain a rough number decays slowly as $n$ increases; see \Cref{gaps-fig}.

\begin{figure}
\centering
\includegraphics[width=0.8\textwidth]{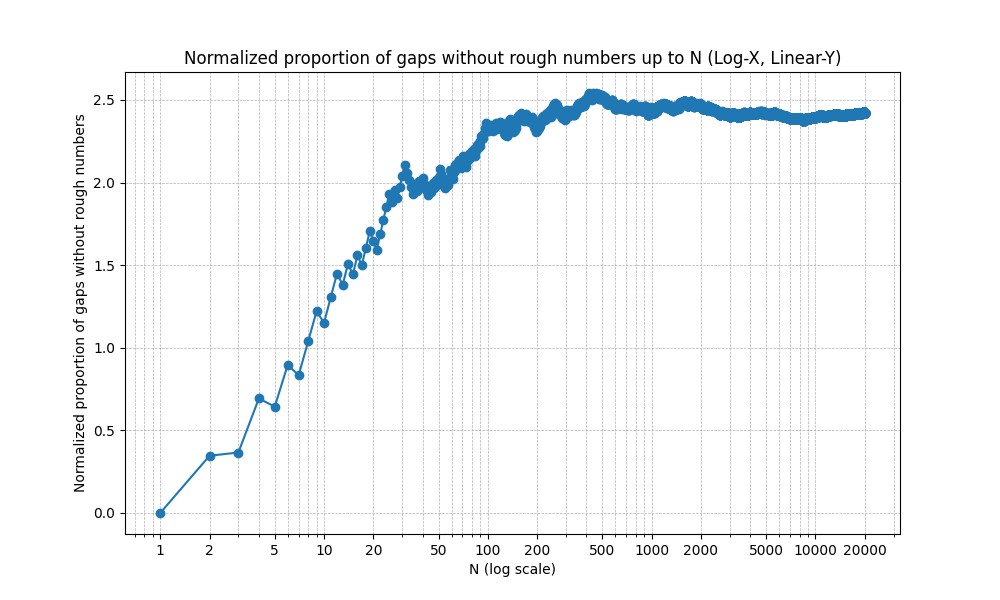}
\caption{The proportions from the previous figure, multiplied by $\log N$.  \Cref{erdos theorem} implies that these products are bounded as $N \to \infty$, and conjecturally converge to a limit $c>0$, which we believe to be between $2.7$ and $2.8$, although this is only weakly supported by the available numerical data.}\label{gaps-norm-fig}
\end{figure}

This set of $n$ was studied by Erd\H{o}s \cite[p. 74]{erdos}, who initially thought that the $n^{\mathrm{th}}$ prime gap would always contain a rough number for sufficiently large $n$, but then wrote:

\begin{quote}
...I am now sure that this is not true and I ``almost'' have a counterexample. Pillai and Szekeres observed that for every $t \leq 16$, a set of $t$ consecutive integers always contains one
which is relatively prime to the others. This is false for $t = 17$, the smallest counterexample being $2184, 2185, \dots, 2200$. Consider now the two arithmetic progressions
$2183 + d \cdot 2 \cdot 3 \cdot 5 \cdot 7 \cdot 11 \cdot 13$ and $2201 + d \cdot 2 \cdot 3 \cdot 5 \cdot 7 \cdot 11 \cdot 13$. There certainly will be infinitely many values of $d$ for which the progressions simultaneously represent primes; this
follows at once from hypothesis H of Schinzel, but cannot at present be proved.  These primes are consecutive and give the required counterexample. I expect that
this situation is rather exceptional and that the integers $k$ for which there is no $m$ satisfying $p_k < m < p_{k+1}$ and $p(m) > p_{k+1} - p_k$ have density $0$.
\end{quote}

In fact, as we shall remark in \Cref{conditional-sec}, one can make the counterexample construction simpler: ``cousin primes'' (consecutive primes differing by $4$) are necessarily of the form $1+6d$ and $5+6d$, and would also produce a prime gap that does not contain a rough number.

There is considerable work in the literature on the distribution of prime numbers in short intervals, as well as the distribution of prime gaps, but a significant portion of it is conditional on unsolved conjectures such as the Hardy--Littlewood prime tuples conjecture.  As we shall see in this paper, though, one can use modern sieve theory to translate some of this work to the setting of rough numbers in short intervals (and specifically, to the intervals associated to prime gaps), and obtain good results even at the unconditional level.  More precisely, in this paper we confirm the latter expectation of Erd\H{o}s (also listed at \cite[Problem \#682]{bloom}):

\begin{theorem}[Erdos \#682]\label{erdos theorem}
For $X>2$, let $N(X)$ be the number of prime gaps $(p_n, p_{n+1})$ with $p_n \in [X,2X]$ that do not contain a rough number (i.e., $p(m) < p_{n+1} - p_n$ for all $p_n < m < p_{n+1}$).  Then one has\footnote{See \Cref{notation-sec} below for our conventions on asymptotic notation.}
\begin{equation}\label{final-bound}
N(X) \ll \frac{X}{\log^2 X}.
\end{equation}
Assuming the Dickson--Hardy--Littlewood prime tuples conjecture in the form \Cref{dhltc} below, we can improve this to
\begin{equation}\label{asymp}
N(X) \sim c \frac{X}{\log^2 X}
\end{equation}
for some (explicitly describable) constant $c>0$.
\end{theorem}

Since the prime number theorem asserts that the number of primes with $X \le p_n \le 2X$ is $\asymp \frac{X}{\log X}$ as $X \to \infty$, we do indeed see (after concatenating together various dyadic intervals $[X,2X]$) that the proportion of the first $N$ prime gaps without rough numbers decays to zero as $N \to \infty$, as predicted by Erd\H{o}s; indeed we obtain an upper bound of $O( 1 / \log N)$ on this rate, and conditionally on the prime tuples conjecture one can sharpen this to $\sim c / \log N$.  It is challenging to compute the constant $c$ to high precision, but we believe it is between $2.7$ and $2.8$, which is (weakly)\footnote{The slow convergence is perhaps unsurprising in view of some lower order irregularities \cite{lemke} in the distribution of prime gaps.} consistent with our preliminary numerics in \Cref{gaps-norm-fig}.

\begin{remark} In the above, Erd\H{o}s stated the problem with the condition \eqref{rough-gap} and also with the slightly stronger condition $p(m) > p_{n+1} - p_n$.  As the gap length $p_{n+1} - p_n$ is prime if and only if it is equal to two, the latter amounts to adding the set of twin primes to the count $N(X)$.  Standard sieve theory shows that the size of this contribution is $O( X / \log^2 X)$, and on the prime tuples conjecture it is $\sim {\mathfrak S} X / \log^2 X$, where the twin primes constant ${\mathfrak S} = 1.3203236\dots$ is defined in \eqref{twinprime}.  Thus, the results of this paper also extend to this slightly modified formulation, after increasing the constant $c$ by the twin prime constant ${\mathfrak S}$.
\end{remark}

\begin{remark} It remains an open problem whether there are infinitely many prime gaps that do not contain rough numbers, or to establish the slightly stronger claim that $N(X)$ is non-zero for all sufficiently large $X$.  This would follow for instance from the Polignac conjecture \cite{polignac}, which in particular implies that there are infinitely many cousin primes.  Unconditionally, our analysis suggests that it will be the prime gaps of very small length that have the best chance of avoiding rough numbers. It is then tempting to apply the machinery of Zhang \cite{zhang} and Maynard \cite{maynard} (see also \cite{polymath}), which can produce infinitely many clusters $n+h_1,\dots,n+h_k$ of almost primes a bounded distance apart (for certain ``admissible'' tuples $h_1,\dots,h_k$), with the numbers between these clusters having small prime factors, and with at least two of the almost primes in the cluster actually being prime (cf. the main theorem of \cite{pintz}).  Unfortunately, the proof (which is based on a weighted pigeonhole principle) does not give control on \emph{which} two of the elements $n+h_1,\dots,n+h_k$ of the cluster are prime, and in particular does not seem to preclude the possibility that any two consecutive primes in the cluster surround one of the almost primes of that cluster, thus failing to provide an example of a prime gap without a rough number.  We therefore suspect that unconditionally resolving the question of whether there are infinitely many prime gaps without rough numbers will require further advances in the theory of small prime gaps.
\end{remark}

One can explain \Cref{erdos theorem} heuristically as follows.  By the prime number theorem, a typical prime gap $(p_n, p_{n+1})$ in $[X,2X]$ would be expected to have length comparable to $\log X$; thus the theorem asserts, roughly speaking, that for $H$ slightly smaller than $\log X$, that a typical interval $(x,x+H]$ in $[X,2X]$ should contain a rough number.  Pick a $z$ somewhat larger than the expected prime gap (e.g., $z = \exp(\log^\beta X)$ for some small constant $\beta$), and consider the numbers $m$ which are rough in the sense that $p(m) \geq z$.  By well-known results of Buchstab \cite{buchstab} and de Bruijn \cite{deBruijn}, the density of such rough numbers is comparable to $\frac{1}{\log z}$; and from standard sieve theory methods, one expects this set of rough numbers to be distributed in a somewhat random (or ``Poissonian'') fashion, after taking into account local irregularities arising from small primes (in the spirit of \cite{gallagher}).  So, if $\log z$ is small compared with $H$, one expects (by the concentration of measure phenomenon) that most intervals $(x,x+H]$ would contain about $H/\log z$ such rough numbers, which should give the claim.

The main technical challenge in making this argument rigorous is to establish a suitable concentration of measure result for rough numbers in short intervals.  The simplest way to proceed is via the second moment method; this partially works, but yields a weaker version of \Cref{erdos theorem} in which the upper bound $N(X) \ll X / \log^2 X$ is worsened to
\begin{equation}\label{weak-bound}
    N(X) \ll X / \log^{4/3-o(1)} X.
\end{equation}
Here one can use well-known bounds on mean value asymptotics of singular series that date back to the work of Montgomery \cite{montgomery}; we give the details in \Cref{firstmoments-sec}.  To recover the near-optimal exponent, one must also control higher moments of rough number counts in short intervals.  This ultimately amounts to controlling averages of $k$-point singular series, a topic first investigated by Gallagher \cite{gallagher} in connection with the closely related question of understanding moments of \emph{prime} number counts in short intervals.  Here it turns out that the more recent calculations of Montgomery and Soundararajan \cite{ms} give a sufficiently strong concentration of measure result\footnote{In fact, their calculations should ultimately give a central limit theorem as well, although we will not perform this calculation here.} to establish a nearly optimal bound
\begin{equation}\label{medium-bound}
    N(X) \ll X / \log^{2-o(1)} X;
\end{equation}
see \Cref{higher-moments-sec}. A refinement of the analysis, that is usable for all but the shortest prime gaps, then recovers the full bound \eqref{final-bound}.  The final asymptotic \eqref{asymp} then follows by applying the prime tuples conjecture to control more precisely the prime gaps of bounded length; see \Cref{conditional-sec} for details.

\subsection*{Acknowledgments}

AG is supported by NSF grant OIA-2229278. TT is supported by NSF grant DMS-2347850, the James and Carol Collins Chair, and the Mathematical Analysis \& Application Research Fund; TT is particularly grateful to recent donors to that fund.  We are indebted to Thomas Bloom for creating and maintaining the Erd\H{o}s problems website \cite{bloom}, where we learned of this problem.

\subsection{Notation}\label{notation-sec}

We use the following asymptotic notation throughout the paper:

\begin{itemize}
\item We use $X \ll Y$, $Y \gg X$, or $X = O(Y)$ to denote the estimate $|X| \leq C Y$ for some constant $C>0$, which can depend on fixed quantities such as the exponents $\alpha,\beta$ appearing in the next two sections, or the parameter $K$ in the final section.  We use $X \asymp Y$ to denote the estimate $X \ll Y \ll X$.
\item We use $Y = o(Z)$ as $X \to \infty$ to denote the estimate $|Y| \leq c(X) Z$ for some function $c(X)$ of a parameter $X$ that goes to zero as $X \to \infty$ (holding all other parameters fixed).  We use $Y \sim Z$ to denote the estimate $Y = (1+o(1)) Z$ as $X \to \infty$.
\end{itemize}

We use $\mathbbm{1}_A$ to denote the indicator function of a set $A$, thus $\mathbbm{1}_A(n) = 1$ if $n \in A$ and $\mathbbm{1}_A(n) = 0$ otherwise.

All sums and products over the variable $p$ are understood to be over primes.

\section{First attempt: using the second moment method}\label{firstmoments-sec}

In this section we establish the weak version \eqref{weak-bound} of \Cref{erdos theorem}, which already suffices to answer the original question of Erd\H{o}s.  This simplified argument is based on the second moment method, and will serve to illustrate the more complicated arguments in the later sections of this paper.

The key input will be the following moment estimates for counting rough numbers in short intervals.

\begin{lemma}[First and second moment estimates]\label{moment-estimates}  Let $0 < \beta < \alpha < 1$ be fixed constants; we permit all implied constants to depend on these parameters.  Let $X \geq 2$, and define the quantities
$$ H \coloneqq \lfloor \log^\alpha X \rfloor; \quad z \coloneqq \exp( \log^\beta X ).$$
Let $\mathcal R(z) \coloneqq \{ n: p(n) \geq z \}$ denote the set of numbers $n$ whose smallest prime factor $p(n)$ is at least $z$, and set
\begin{equation}\label{M-def}
    R \coloneqq \frac{e^{-\gamma} H}{\log z} \asymp \log^{\alpha-\beta} X,
\end{equation}
where $\gamma$ is the Euler--Mascheroni constant.  Then one has the first moment estimate
\begin{equation}\label{first-moment}
\frac{1}{X} \int_X^{2X} \sum_{x < n \le x+H} \mathbbm{1}_{\mathcal R(z)}(n) \, dx = R + O\left(\frac{1}{\log^2 X}\right)
\end{equation}
and the second moment estimate
\begin{equation}\label{second-moment}
\frac{1}{X} \int_X^{2X} \left|\sum_{x < n \le x+H} \mathbbm{1}_{\mathcal R(z)}(n)\right|^2 \, dx = R^2 + O(R).
\end{equation}
In particular, we have
\begin{equation}\label{variance}
 \frac{1}{X} \int_X^{2X} \left|\sum_{x < n \le x+H} \mathbbm{1}_{\mathcal R(z)}(n) - R \right|^2 \, dx \ll R.
\end{equation}
\end{lemma}

Informally, \eqref{variance} asserts that most intervals $[x, x+H]$ with $x \in [X,2X]$ contain approximately $R$ elements of $\mathcal R(z)$.

\begin{proof}  We begin with the first moment estimate \eqref{first-moment}. We may assume $X$ to be larger than any given constant (depending on $\alpha,\beta$). We have
\begin{align*}
 \frac{1}{X} \int_X^{2X} \sum_{x< n \le x+H} \mathbbm{1}_{\mathcal R(z)}(n)  \, dx & =  \frac{1}{X}\sum_{X\le n \le 2X+H} \mathbbm{1}_{\mathcal R(z)}(n)  \int_{[X,2X]\cap [n-H, n)} 1 \, dx \\&
= \frac{H}{X}\sum_{X\le n \le 2X} \mathbbm{1}_{\mathcal R(z)}(n)  + O\left(\frac{H^2}{X}\right).
\end{align*}
From Buchstab's theorem \cite[Theorem 7.11]{mont-vaughan} we have
$$ \sum_{n < X} \mathbbm{1}_{\mathcal R(z)}(n) = \frac{\omega(u) X}{\log z} - \frac{z}{\log z} + O\left( \frac{X}{\log^2 X} \right)$$
where $u \coloneqq \log X / \log z = \log^{1-\beta} X$ and $\omega$ is the Buchstab function.  Using the de Bruijn asymptotics \cite{deBruijn}
$$ \omega(u) = e^{-\gamma} + O\left( \frac{1}{\Gamma(u+1)} \right),$$
where $\Gamma$ denotes the Gamma function, one then easily concludes that
$$ \sum_{n < X} \mathbbm{1}_{\mathcal R(z)}(n) = \frac{e^{-\gamma} X}{\log z} + O\left( \frac{X}{\log^2 X} \right)$$
and a similar computation (with a very slightly different value of $u$) gives
$$ \sum_{n \leq 2X} \mathbbm{1}_{\mathcal R(z)}(n) = \frac{2e^{-\gamma} X}{\log z} + O\left( \frac{X}{\log^2 X} \right)$$
and hence
$$\sum_{X\le n \le 2X} \mathbbm{1}_{\mathcal R(z)}(n) = \frac{e^{-\gamma} X}{\log z} + O\left( \frac{X}{\log^2 X} \right).$$
Inserting this back into the previous estimates, we obtain \eqref{first-moment} after a brief calculation.

Now we turn to the second moment estimate \eqref{second-moment}.  We expand
\begin{align*}
& \frac 1X \int_X^{2X}\left(\sum_{x\le n \le x+H} \mathbbm{1}_{\mathcal R(z)}(n)\right)^2 \, dx  \\
& =  \frac 1X  \int_X^{2X} \sum_{x\le n \le x+H}  \sum_{x\le m \le x+H} \mathbbm{1}_{\mathcal R(z)}(n)\mathbbm{1}_{\mathcal R(z)}(m)\, dx  \\
 & = \frac 1X \int_X^{2X} \sum_{x\le n \le x+H}  \mathbbm{1}_{\mathcal R(z)}(n) \, dx
+ \frac 2X \int_X^{2X} \sum_{1\le h \le H}  \sum_{x\le n \le x+H-h} \mathbbm{1}_{\mathcal R(z)}(n)\mathbbm{1}_{\mathcal R(z)}(n+h) \\
& = O(R) + \frac 2X  \sum_{1\le h \le H}  \sum_{X\le n \le 2X+H-h} \mathbbm{1}_{\mathcal R(z)}(n)\mathbbm{1}_{\mathcal R(z)}(n+h) \int_{[X,2X]\cap [n-H+h, n]} 1 \, dx  \\
& = \frac 2X  \sum_{1\le h \le H} (H-h)  \sum_{X\le n \le 2X} \mathbbm{1}_{\mathcal R(z)}(n)\mathbbm{1}_{\mathcal R(z)}(n+h) + O(R)
\end{align*}
where we have used \eqref{first-moment} to control the diagonal terms; some error terms are generated at the boundary $n = X + O(H)$ or $n = 2X + O(H)$ of the interval $[X,2X]$, but their contribution easily seen to be negligible compared to the existing error term $O(R)$ (they are of size $O(\frac{H^2}{X})$, which is far smaller).  Setting
$$S_h \coloneqq \frac{1}{X} \sum_{X\le n \le 2X} \mathbbm{1}_{\mathcal R(z)}(n)\mathbbm{1}_{\mathcal R(z)}(n+h),$$
we can then rewrite the left-hand side of \eqref{second-moment} as
$$
2\sum_{1\le h \le H} (H-h) S_h + O(R).
$$
As all the elements of $\mathcal{R}(z)$ are odd, we can restrict $h$ to be even, and rewrite the above expression as
\begin{equation}\label{moment-breakdown}
2\sum_{1\le m \le H/2} (H-2m) S_{2m} + O(R).
\end{equation}

Let $1 \leq m \leq H/2$.   We will estimate $S_{2m}$ using the Rosser--Iwaniec sieve. Set $P(z) \coloneqq \prod_{p<z} p$.   Define a sequence $a_{2m} \colon \N \to \R$ by setting $a_{2m}(\ell) \coloneqq 1$ if $\ell = n(n+2m)$ for some $X \le n \le 2X$, and $a_{2m}(\ell) \coloneqq 0$ otherwise.  Then we have
 \begin{equation*}
 S_{2m} = \frac{1}{X} \sum_{\substack{\ell\\ (\ell, P(z))=1}} a_{2m}(\ell)
 \end{equation*}
 where $(\ell,P(z))$ denotes the greatest common divisor of $\ell$ and $P(z)$.
From the Chinese remainder theorem we can write
 \begin{equation*}
 \sum_{\ell}  a_{2m}(d\ell) = \frac{X\rho_{2m}(d)}{d} + O(2^{\omega(d)}),
 \end{equation*}
 for any squarefree $d$, where $\omega(d)$ denotes the number of prime divisors of $d$, and $\rho_{2m}$ is a multiplicative function with
 \begin{equation*}
 \rho_{2m}(p) = \begin{cases}
 1-\frac2p, & p\nmid 2m,\\
  1-\frac1p, & p\mid 2m.
 \end{cases}
 \end{equation*}
 Define the quantity $V_{2m}(z)$ by the formula
 \begin{equation*}
V_{2m}(z) \coloneqq \prod_{p<z} \rho_{2m}(p) = \frac12   V(z) \prod_{\substack{ p>2 \\ p\mid m }} \left(\frac{p-1}{p-2}\right)
\end{equation*}
where
$$ V(z) \coloneqq \prod_{2<p<z} \left(1-\frac{2}{p}\right).$$

To apply the sieve estimate, we need a bound on $\frac{V_{2m}(w)}{V_{2m}(z)}$ for $w\le z$.  Here we will use a somewhat crude estimate:
\begin{align*}
\frac{V_{2m}(w)}{V_{2m}(z)} & = \prod_{w \leq p< z} \rho_{2m}(p)\\
& \le  \prod_{w \leq p< z} \left(1-\frac{2}{p}\right)^{-1} \prod_{ p \le 2m} \left(1-\frac{1}{p}\right)^{-1}  \\
& \ll \left(\frac{\log z}{\log w}\right)^2 (\log 2m + O(1)) \\
& \ll (\log^{o(1)} X) \left(\frac{\log z}{\log w}\right)^2
\end{align*}
thanks to Mertens' theorem \cite[Theorem 2.7]{mont-vaughan}.

We are now ready to apply the sieve.  Let $s\ge 19$ be a parameter to be chosen later.  By \cite[Theorem 6.9]{fried-iwa} (and crudely bounding $2^{\omega(m)}$ by $m$), we have
\begin{align*}
S_{2m} & = V_{2m}(z)(1+ E(X,z,m,s))  + O\bigg(\frac{1}{X} \sum_{\substack{m\le z^s \\ m|P(z)}} 2^{\omega(m)} \bigg) \nonumber\\
& = V_{2m}(z)(1+ E(X,z,m,s)) + O\left( \frac{z^{2s}}{X} \right)
\end{align*}
where the first error term $E(X,z,m,s)$ obeys a bound of the form
\begin{equation}\label{error-bound}
 |E(X,z,m,s)| \ll e^{-s} \log^{o(1)} X.
\end{equation}
If we then choose $s = \log^{(1-\beta)/2} X$ (say), both error terms can be shown to be $O(\frac{1}{\log^{10} X})$ (with substantial room to spare), thus
$$ S_{2m} = V_{2m}(z) + O\left( \frac{1}{\log^{10} X} \right).$$
The quantity \eqref{moment-breakdown} then simplifies to
\begin{equation}\label{moment-breakdown-2}
V(z) \sum_{1 \leq m \leq H/2} (H-2m) \prod_{\substack{ p>2 \\ p\mid m }} \left(\frac{p-1}{p-2}\right) + O(R).
\end{equation}

From Mertens' theorem \cite[Theorem 2.7(e)]{mont-vaughan} one has the approximation
\begin{equation}\label{Merten-0}
 V(z)  = 2 \frac{\mathfrak{S}e^{-2\gamma}}{\log^2 z} \left(1+ O\left(\frac{1}{\log z}\right) \right),
\end{equation}
where $\mathfrak{S}$ is the twin prime constant
\begin{equation}\label{twinprime}
\mathfrak{S} \coloneqq 2 \prod_{p>2}\left(1-\frac{1}{p}\right)^{-2} \left(1-\frac{2}{p}\right) =  2 \prod_{p>2}\left(1 - \frac{1}{(p-1)^2} \right) =1.3203236 \dots
\end{equation}
Actually, the error term in \eqref{Merten-0} can be improved using the prime number theorem, and specifically we have
\begin{equation}\label{Merten-1}
 V(z)  = \frac{2\mathfrak{S}e^{-2\gamma}}{\log^2 z} \left(1+ O\left(\exp( - c \sqrt{\log z} ) \right) \right)
\end{equation}
for some absolute constant $c>0$.  This type of improvement is well-known folklore (and one can do even better by using the Vinogradov--Korobov error term in the prime number theorem, or by assuming the Riemann hypothesis), but for the convenience of the reader we include a proof here.  From the prime number theorem with classical error term we have
$$ \sum_{w < p \leq z} \frac{1}{p} = \log\log z - \log\log w + O(\exp(-c \sqrt{\log w}))$$
for all $100 \leq w \leq z$ and some $c>0$ (see e.g, \cite[Exercise 6.2.3(a)]{mont-vaughan}); exponentiating this, we soon find that
$$ \frac{V(w)}{V(z)} = \frac{\log^2 z}{\log^2 w} \left( 1 + O(\exp(-c \sqrt{\log w})) \right).$$
Multiplying by $V(z)$, using \eqref{Merten-0}, and then taking the limit as $z \to \infty$, we recover \eqref{Merten-1} (after relabeling $w$ as $z$).  With this improved error term, we can now simplify \eqref{moment-breakdown-2} further as
\begin{equation}\label{moment-breakdown-3}
\frac{2\mathfrak{S}e^{-2\gamma}}{\log^2 z} \sum_{1 \leq m \leq H/2} (H-2m) \prod_{\substack{ p>2 \\ p\mid m }} \left(\frac{p-1}{p-2}\right) + O(M).
\end{equation}
The sum appearing here has been extensively studied in the literature \cite{croft}, \cite{goldston}, \cite{sg}, \cite{montgomery}, \cite{ms-beyond}, \cite{ms}, \cite{vaughan}.  For our application it will suffice to apply the estimates from \cite[(12), (13)]{montgomery}, which assert that\footnote{Sharper estimates can be found in the more recent literature cited above, but this does not lead to a significantly improved bound in our application, as this contribution to the error term is already of lower order.}
\begin{align}
\sum_{m\le w}  \prod_{\substack{ p>2 \\ p\mid m }} \left(\frac{p-1}{p-2}\right) & = \frac{2w}{\mathfrak{S}} + O(\log w) \label{singular-series}\\
\sum_{m\le w} m \prod_{\substack{ p>2 \\ p\mid m }} \left(\frac{p-1}{p-2}\right) & = \frac{w^2}{\mathfrak{S}} + O(w \log w)\nonumber
\end{align}
for any $w > 2$.  A brief calculation then gives
$$
\sum_{1\le m \le \frac{H}{2}} (H-2m)  \prod_{\substack{ p>2 \\ p\mid m }} \left(\frac{p-1}{p-2}\right) = \frac{H^2}{2\mathfrak{S}} + O(H \log H),
$$
and so the expression in \eqref{moment-breakdown-3} simplifies to
$$ R^2 + O(R) + O\left( \frac{H \log H}{\log^2 z} \right) = R^2 + O(R),$$
giving \eqref{second-moment} as required.

Finally, \eqref{variance} follows from \eqref{second-moment} and \eqref{first-moment} after expanding out the square.
\end{proof}

We are now ready to prove the weak version of \Cref{erdos theorem}.  Let $X$ be large, and let $\alpha$, $\beta$, $H$, $z$ be as in \Cref{moment-estimates}.  We define the exceptional set $E(X) \subset [X,2X]$ to be the set of all $x \in [X,2X]$ such that the interval $(x,x+H]$ contains no numbers in $\mathcal R(z)$.  By \Cref{moment-estimates}, we can upper bound the Lebesgue measure $|E(X)|$ of this set by
\begin{align}
|E(X)| & = \frac{1}{R^2} \int_{E(X)}   \left| \sum_{x < n \le x+H} \mathbbm{1}_{\mathcal R(z)}(n) - R\right|^2 \, dx \nonumber \\
 & \le  \frac{1}{R^2} \int_X^{2X} \left| \sum_{x < n \le x+H} \mathbbm{1}_{\mathcal R(z)}(n) - R\right|^2 \, dx \nonumber\\
 &\ll \frac{X}{R}. \label{E-bound}
\end{align}
We now make the following key observation: if $X \leq p_n \leq 2X$ is such that the prime gap $(p_n, p_{n+1})$ does \emph{not} contain a rough number, then at least one of the following statements must hold:
\begin{itemize}
    \item[(i)] (Small gap) One has $p_{n+1} - p_n \leq 2H$.
    \item[(ii)] (Large gap) One has $p_{n+1} - p_n \geq z$.
    \item[(iii)] (Exceptional)  The intersection of $E(z)$ with $(p_n, p_{n+1})$ has measure at least $H$.
    \item[(iv)] (Boundary) One has $p_{n+1} > 2X$.
\end{itemize}

Indeed, suppose that (i), (ii), (iv) all failed, thus $X \leq p_n < p_{n+1}-2H \leq p_{n+1} \leq 2X$ and $p_{n+1} - p_n < z$.  Since the prime gap $(p_n, p_{n+1})$ does not contain a rough number, it must be disjoint from $\mathcal{R}(z)$ since $z$ exceeds the length of the gap. By definition of $E(X)$, we conclude that the interval $[p_n, p_{n+1}-H)$ is contained in $E(X)$.  As this interval has measure at least $H$, we obtain the conclusion (iii).

We can now bound the number $N(X)$ of exceptional $n$ in \Cref{erdos theorem} by considering the scenarios (i)-(iv) separately.  The scenario (iv) contributes at most one case to $N(X)$.  From \eqref{E-bound} we see that scenario (iii) can contribute at most
$$
\ll \frac{1}{H} \frac{X}{R} \asymp \frac{X}{\log^{2\alpha-\beta} X}$$
cases to $N(X)$.  Since the sum of all the prime gaps $p_{n+1}-p_n$ excluding the boundary case (iv) telescopes to at most $X$, the scenario (ii) contributes at most
$$ \ll \frac{X}{z}$$
new cases to $N(X)$.  Finally, to control the contribution of scenario (i) to $N(X)$, we use sieve theory.  For any fixed gap $1 \leq h \leq 2H$, we know from standard sieve-theoretic estimates (see e.g., \cite[Corollary 3.14]{mont-vaughan}) that the number of primes $p_n$ with $X \leq p_n \leq 2X$ and $p_{n+1} - p_n = h$ is bounded by
\begin{equation}\label{nbound}
    O\left( \frac{X}{\log^2 X} \prod_{\substack{ p>2 \\ p\mid h }} \left(\frac{p-1}{p-2}\right)\right).
\end{equation}
Summing this in $h$ using \eqref{singular-series}, we conclude that the total number of contributions of case (i) to $N(X)$ is at most
$$ \ll \frac{X}{\log^2 X} H \ll \frac{X}{\log^{2-\alpha} X}.$$
We conclude that for any fixed $0 < \beta < \alpha < 1$, we have
$$ N(X) \ll \frac{X}{\log^{2-\alpha} X} + \frac{X}{z} + \frac{X}{\log^{2\alpha-\beta} X} + 1.$$
Taking $\alpha = \frac{2}{3}$ and $0 < \beta < \alpha$, we have thus obtained the bound $N(X) \ll X / \log^{4/3-\beta} X$.  Since $\beta > 0$ can be taken arbitrarily small, we obtain the claim \eqref{weak-bound}.

\begin{remark} A variant of the argument (choosing $\alpha$ and $\beta$ both close to $1$) shows that Erd\H{o}s's prediction is in fact true in a rather strong sense: for almost all $n$, the prime gap $(p_n, p_{n+1})$ contains a natural number $m$ with $p(m)$ exceeding $\exp( (p_{n+1}-p_n)^{1-o(1)} )$ (not just exceeding $p_{n+1}-p_n$).  We leave the details of this variant to the interested reader.
\end{remark}

\section{Second attempt: using higher moments}\label{higher-moments-sec}

We now improve \eqref{weak-bound} to \eqref{medium-bound}.  The main idea here is to control higher moments beyond the first and second moments. Specifically, we will need the following variant of \Cref{moment-estimates}, inspired by \cite[Theorem 3]{ms}:

\begin{lemma}[Higher moment estimates]\label{higher-estimates}  Let $\beta$, $\alpha$, $X$, $H$, $z$, $\mathcal{R}(z)$, $R$ be as in \Cref{moment-estimates}.  Let $K$ be a fixed natural number; we allow implied constants to depend on $K$.  Then
\begin{equation}\label{variance-k}
 \frac{1}{X} \int_X^{2X} \left(\sum_{x< n \le x+H} \mathbbm{1}_{\mathcal R(z)}(n) - R' \right)^{K} \, dx \ll H^{K/2+o(1)}
\end{equation}
where $R'$ is the slight modification of $R$ defined by
$$ R' \coloneqq H \prod_{p<z} \left(1 - \frac{1}{p} \right) = R ( 1 + O( \exp(-c\sqrt{\log z})))$$
with the last estimate coming from the improved error term in Mertens' theorem (cf. \eqref{Merten-1}).
\end{lemma}

Let us assume this lemma for the moment and see how it establishes \eqref{medium-bound}.  Let $K$ be a large but fixed even integer.  We repeat the previous arguments, controlling the number $N(X)$ of exceptions to \Cref{erdos theorem} by the four scenarios (i)-(iv) as before.  The only distinction here is that we can now improve \eqref{E-bound} to
\begin{align*}
|E(X)| & = \frac{1}{(R')^K} \int_{E(X)}   \left| \sum_{x < n \le x+H} \mathbbm{1}_{\mathcal R(z)}(n) - R'\right|^K \, dx \\
 & \le  \frac{1}{(R')^K} \int_X^{2X} \left| \sum_{x < n \le x+H} \mathbbm{1}_{\mathcal R(z)}(n) - R'\right|^K \, dx \\
 &\ll \frac{X H^{K/2+o(1)}}{R^K}.
\end{align*}
Tracing through the argument, we conclude that the number of $n$ contributing to scenario (iii) now has the improved bound
$$
\ll \frac{1}{H} \frac{X H^{K/2+o(1)}}{R^{K}} \asymp \frac{X}{\log^{(K+2)\alpha/2-K\beta-o(1)} X}$$
and so now the total bound on exceptional $n$ is
\begin{equation}\label{nx}
    N(X) \ll \frac{X}{\log^{2-\alpha} X} + \frac{X}{z} + \frac{X}{\log^{(K+2)\alpha/2-K\beta-o(1)} X} + 1.
\end{equation}
If we set $\alpha = \frac{4}{K+4}$ and $\beta = \frac{1}{K^2}$, we have thus bounded the number $N(X)$ of exceptions by $O( X / \log^{2 - \frac{4}{K+4} - \frac{1}{K} - o(1)} X)$.  Taking $K$ arbitrarily large, we obtain the desired bound of $O( X / \log^{2-o(1)} X)$.

It remains to establish \Cref{higher-estimates}.  The cases $K=0,1,2$ are already covered by \Cref{moment-estimates} (with slightly better bounds), so we may take $K > 2$.  We will follow the combinatorial calculations from \cite[Section 3]{ms}, combined with the singular series estimates from that paper as well as the Rosser--Iwaniec sieve.  Our analysis will be slightly simpler than the corresponding one in \cite{ms} because the rough numbers $\mathcal R(z)$ have a constant density $\frac{R'}{H}$ throughout the interval $[X,2X]$, whereas the primes in $[X,2X]$ have an (approximate) density function $\frac{1}{\log x}$ that slowly varies in $[X,2X]$; we can therefore avoid some technical manipulations in \cite{ms} designed to deal with this variable density.

We introduce the von Mangoldt-like function
\begin{equation}\label{vm-def}
     \Lambda^{(z)}(n) \coloneqq \frac{H}{R'} \mathbbm{1}_{\mathcal R(z)}(n)
\end{equation}
and its normalized variant
$$ \Lambda^{(z)}_0(n) \coloneqq \Lambda^{(z)}(n) - 1.$$
These functions will play the role of $\Lambda$ and $\Lambda_0$ in \cite[Section 3]{ms}.  We can renormalize \eqref{variance-k} as
$$ \int_X^{2X} \left(\sum_{x< n \le x+H} \Lambda^{(z)}_0(n) \right)^{K} \, dx \ll \left(\frac{H}{R'}\right)^K H^{K/2+o(1)} X.$$
The integrand is piecewise constant on unit intervals $[n',n'+1)$, so it suffices to show that
\begin{equation}\label{hr}
     \sum_{X \leq n' \leq 2X} \left(\sum_{n'< n \le n'+H} \Lambda^{(z)}_0(n) \right)^{K} \ll \left(\frac{H}{R'}\right)^K H^{K/2+o(1)} X
\end{equation}
since the boundary contributions from $x = X+O(1)$ or $x = 2X+O(1)$ are negligible.  Expanding out the $K^{\mathrm{th}}$ power, we can rewrite the left-hand side as
$$ \sum_{1 \leq d_1,\dots,d_k \leq H} \sum_{X \leq n' \leq 2X} \Lambda^{(z)}_0(n'+d_1) \dots \Lambda^{(z)}_0(n'+d_k).$$
The $d_1,\dots,d_k$ are \emph{not} currently required to be distinct, which forces us to perform some combinatorial rearranging. Following \cite[(61)]{ms}, we regroup the above expression (replacing $n'$ with $n$) as
\begin{equation}\label{bigmess}
     \sum_{k=1}^K \sum_{\stackrel{M_1,\dots,M_k \geq 1}{\sum M_i = K}} \binom{K}{M_1 \dots M_k} \frac{1}{k!} \sum_{\stackrel{1 \leq d_1,\dots,d_k \leq H}{d_i \text{ distinct}}} \sum_{X \leq n \leq 2X} \prod_{i=1}^k \Lambda^{(z)}_0(n+d_i)^{M_i}
\end{equation}
where the $M_i$ range over positive integers summing to $H$ and $\binom{K}{M_1 \dots M_k}  \coloneqq \frac{K!}{M_1! \dots M_k!}$.  Introducing the functions
$$ \Lambda^{(z)}_m(n) \coloneqq \Lambda^{(z)}(n)^m \Lambda^{(z)}_0(n) $$
we can repeat the derivation of \cite[(62), (63)]{ms} to expand \eqref{bigmess} further as
\begin{align*}
    &\sum_{k=1}^K \frac{1}{k!} \sum_{\stackrel{M_1,\dots,M_k \geq 1}{\sum M_i = K}} \binom{K}{M_1 \dots M_k}  \\
    &\quad \times \sum_{\stackrel{m_1,\dots,m_k}{0 \leq m_i < M_i}} \prod_{i=1}^k (-1)^{M_i-1-m_i} \binom{M_i-1}{m_i} L_k(\mathbf{m})
\end{align*}
where
\begin{equation}\label{lk}
 L_k(\mathbf{m}) \coloneqq \sum_{\stackrel{1 \leq d_1,\dots,d_k \leq H}{d_i \text{ distinct}}} \sum_{X \leq n \leq 2X} \prod_{i=1}^k \Lambda^{(z)}_{m_i}(n+d_i).
\end{equation}
Since we allow implicit constants to depend on $K$, the number of possible $k$, $M_i$, $m_i$, as well as the multinomial coefficients $\binom{K}{M_1 \dots M_k}$, are bounded by $O(1)$, so by the triangle inequality it will suffice to show that
\begin{equation}\label{hr-lk}
L_k(\mathbf{m}) \ll \left(\frac{H}{R'}\right)^K H^{K/2+o(1)} X
\end{equation}
for all $\mathbf{m} = (m_1,\dots,m_k)$ appearing in the above expansions.

Fix $\mathbf{m}$.  Following \cite{ms}, we introduce the sets
$${\mathcal K} \coloneqq \{1,\dots,k\}; \quad {\mathcal I} \coloneqq \{ i \in {\mathcal K}: m_i = 0 \}$$
and observe that
$$ \prod_{i=1}^k \Lambda^{(z)}_{m_i}(n+d_i) = \left( \frac{H}{R'} \right)^{\sum_i m_i} \sum_{{\mathcal I} \subset {\mathcal J} \subset {\mathcal K}} \prod_{i \in {\mathcal J}} \Lambda^{(z)}_0(n+d_i).$$
It thus becomes relevant to study the sums
\begin{equation}\label{norm}
    \sum_{X \leq n \leq 2X} \prod_{i \in {\mathcal J}} \Lambda^{(z)}_0(n+d_i).
\end{equation}
We first consider the unnormalized sum
\begin{equation}\label{unnorm}
    \sum_{X \leq n \leq 2X} \prod_{i \in {\mathcal J}} \Lambda^{(z)}(n+d_i).
\end{equation}
From \eqref{vm-def}, this is
$$ \left( \frac{H}{R'} \right)^{|{\mathcal J}|} \sum_{X \leq n \leq 2X} \prod_{i \in {\mathcal J}} \mathbbm{1}_{\mathcal R(z)}(n+d_i)$$
where $|{\mathcal J}|$ denotes the cardinality of ${\mathcal J}$.
Invoking the Rosser--Iwaniec sieve \cite[Theorem 6.9]{fried-iwa} as in the previous section, we can express this as\footnote{Recall that the implicit constants in our asymptotic notation are permitted to depend on $K$, and hence on quantities bounded by $K$ such as $|{\mathcal J}|$.}
$$ \left( \frac{H}{R'} \right)^{|{\mathcal J}|} \left( X {\mathfrak S}^{<z}({\mathcal D}_{\mathcal J}) \prod_{p \leq z} \left(1 - \frac{1}{p}\right)^{|{\mathcal J}|}
(1+ E(X,z,{\mathcal D}_{\mathcal J},s)) + O( z^{Ks} )\right)$$
for any $s \geq 19$, where the error term $E(X,z,{\mathcal D}_{\mathcal J},s)$ obeys a bound of the form \eqref{error-bound},and
$$ {\mathfrak S}^{<z}({\mathcal D}_{\mathcal J}) \coloneqq \prod_{p<z} \left( 1-\frac{1}{p} \right)^{-|{\mathcal J}|} \left( 1 - \frac{\nu_p({\mathcal D}_{\mathcal J})}{p} \right)$$
is the truncated version of the singular series
$$ {\mathfrak S}({\mathcal D}_{\mathcal J}) \coloneqq \prod_p \left( 1-\frac{1}{p} \right)^{-|{\mathcal J}|} \left( 1 - \frac{\nu_p({\mathcal D}_{\mathcal J})}{p} \right)$$
and where for each prime $p$, $\nu_p({\mathcal D}_{\mathcal J})$ is the number of distinct residue classes modulo $p$ found amongst ${\mathcal D}_{\mathcal J} \coloneqq \{d_j: j \in {\mathcal J}\}$.  For $p>z$, one has $\nu_p({\mathcal D}_{\mathcal J}) = |{\mathcal J}|$, so one has a good agreement between the singular series and its truncation:
\begin{equation}\label{sdj}
 {\mathfrak S}^{<z}({\mathcal D}_{\mathcal J}) =  {\mathfrak S}({\mathcal D}_{\mathcal J}) \left( 1 + O\left(\frac{1}{z}\right) \right).
\end{equation}
So if we select $s = \log^{(1-\beta)/2} X$ as before, we can evaluate the unnormalized sum \eqref{unnorm} after a brief calculation as
$$ X {\mathfrak S}^{<z}({\mathcal D}_{\mathcal J}) + O\left( \frac{X}{\log^{10 K} X} \right) = X {\mathfrak S}({\mathcal D}_{\mathcal J}) + O\left( \frac{X}{\log^{10 K} X} \right)$$
(say); this is an (unconditional) analogue of the Hardy--Littlewood prime tuples conjecture for rough numbers.  Writing $\Lambda^{(z)}_0 = \Lambda^{(z)}-1$, we can then evaluate the normalized sum \eqref{norm} as
$$ X {\mathfrak S}_0({\mathcal D}_{\mathcal J}) + O\left( \frac{X}{\log^{10 K} X} \right)$$
where, following \cite{ms}, we define
$$ {\mathfrak S}_0({\mathcal D}_{\mathcal J}) \coloneqq \sum_{{\mathcal J'} \subset {\mathcal J}} (-1)^{|{\mathcal J'}|} {\mathfrak S}({\mathcal D}_{\mathcal J'})$$
with the convention ${\mathfrak S}(\emptyset) = 1$.  Inserting this into \eqref{lk}, we conclude that
$$ L_k(\mathbf{m}) \ll X \left| \sum_{\stackrel{1 \leq d_1,\dots,d_k \leq H}{d_i \text{ distinct}}} {\mathfrak S}_0({\mathcal D}_{\mathcal J}) \right| + \frac{X}{\log^{9 K} X}.$$
Applying \cite[Theorem 2]{ms} to sum over the indices $d_j$ with $j \in {\mathcal J}$, and noting that the number of possibilities for the remaining indices is independent of those $d_j$, and comparable to $H^{k - |{\mathcal J}|}$, we conclude that
$$ \sum_{\stackrel{1 \leq d_1,\dots,d_k \leq H}{d_i \text{ distinct}}} {\mathfrak S}_0({\mathcal D}_{\mathcal J}) \ll H^{k - |{\mathcal J}|} H^{|{\mathcal J}|/2+o(1)} = H^{k - |{\mathcal J}|/2 + o(1)}$$
(in fact significantly sharper asymptotics than this are obtained, but the upper bound suffices for our application).  Inserting this back into \eqref{lk}, we conclude that
$$ L_k(\mathbf{m}) \ll X \left( \frac{H}{R} \right)^{\sum_i m_i} \sum_{{\mathcal I} \subset {\mathcal J} \subset {\mathcal K}} H^{k - |{\mathcal J}|/2 + o(1)}.$$
The set ${\mathcal I}$ omits at most $\sum_i m_i$ values from ${\mathcal K}$, thus
$$ |{\mathcal J}| \geq |{\mathcal I}| \geq k - \sum_i m_i.$$
Since $m_i \leq M_i-1$ and $\sum_{i=1}^k M_i = K$, we also have $\sum_i m_i \leq K - k$, thus
$$ |{\mathcal J}| \geq 2k - K.$$
We conclude that
$$ L_k(\mathbf{m}) \ll X \left( \frac{H}{R} \right)^{K-k} H^{K/2  + o(1)} \ll \left( \frac{H}{R} \right)^{K} H^{K/2  + o(1)},$$
thus giving \eqref{hr-lk} as claimed. This concludes the proof of \Cref{higher-estimates}, and hence \eqref{medium-bound}.

\section{Refined analysis}\label{conditional-sec}

We can push the above analysis further and obtain the remaining portions of \Cref{erdos theorem}.  If we remove the contribution $O(X/\log^{2-\alpha} X)$ of the case (i) terms to \eqref{nx}, and select $K$ and $\beta$ appropriately, we conclude that for any fixed $\alpha > 0$, the number of prime gaps $(p_n, p_{n+1})$ of length at least $\log^\alpha X$ with $p_n \in [X,2X]$ that do not contain a rough number is bounded by $O( X / \log^A X)$ for any fixed $A>0$.

For a given gap length $h \geq 2$, let $N_h(X)$ denote the number of prime gaps $(p_n, p_{n+1})$ of length exactly $h$ with $p_n \in [X,2X]$ that do not contain a rough number.  By the above discussion we have
$$ N(X) = \sum_{2 \leq h \leq \log^\alpha X} N_h(X) + O\left(\frac{X}{\log^A x} \right)$$
for any fixed $\alpha,A>0$.  In particular we have
\begin{equation}\label{no}
N(X) = \sum_{2 \leq h \leq \log^{1/100} X} N_h(X) + o\left(\frac{X}{\log^2 x} \right).
\end{equation}
It is now of interest to study the quantities $N_h(X)$.  We can modify the previous analysis to show that $N_h(X)$ exhibits good decay in $h$:

\begin{lemma}\label{Nh-bound}  For $2 \leq h \leq \log^{1/100} X$, we have $N_h(X) \ll \frac{X}{h^2 \log^2 X}$.
\end{lemma}

In fact we can obtain any negative power of $h$ on the right-hand side, but the quadratic decay $1/h^2$ is already summable in $h$ and thus suitable for our purposes.

\begin{proof}  We may assume $X$ to be large. Clearly $N_h(X)$ vanishes for odd $h$, so we may also assume $h$ to be even.  We set
$$ H \coloneqq h-1; \quad z \coloneqq H^{10}$$
and let ${\mathcal R}(z)$ be as before. If a prime gap $(p_n, p_{n+1})$ contributes to $N_h(X)$, then $p_{n+1} = p_n + h$, and for every $i=1,\dots,H$, $p_n+i$ must be divisible by some prime less than $h$. In particular, recalling that $P(z) \coloneqq \prod_{p<z} p$, we see that $p_n$ must belong to a residue class $p_n \equiv b \hbox{ mod } P(z)$ with
$$ \sum_{b< n \le b+H} \mathbbm{1}_{\mathcal R(z)}(n) = 0.$$
For each such residue class, the number of primes $p \in [X, 2X]$ with $p+h$ also prime, and $p = b \hbox{ mod } P(z)$, can be crudely upper bounded by standard sieve estimates (e.g., \cite[Corollary 3.14]{mont-vaughan}) to be
$$ O\left( H^{o(1)} \frac{X}{P(z) \log^2 X} \right)$$
using Mertens' theorem (and $z=H^{O(1)}$) to control the singular series, where the $o(1)$ decay is now with respect to the limit $h \to \infty$.  Since $z = H^{10} \leq \log^{1/10} X$, $P(z)$ is much smaller than $X$, and the previous bound is $O( H^{o(1)} / \log^2 X)$ of the size of the entire residue class $p = b \hbox{ mod } P(z)$ within $[X,2X]$.  We conclude that
$$ N_h(X) \ll \frac{H^{o(1)}}{\log^2 X} | \{ X \leq n' \leq 2X: \sum_{n'< n \le n'+H} \mathbbm{1}_{\mathcal R(z)}(n) = 0 \} |.$$
In particular, we can control $N_h$ in terms of a sixth moment, in that
$$ N_h(X) \ll \frac{1}{H^{6-o(1)} \log^2 X} \int_X^{2X} \left( \sum_{x< n \le x+H} \mathbbm{1}_{\mathcal R(z)}(n) - R' \right)^6\ dx$$
where
$$ R' \coloneqq H \prod_{p<z} \left(1 - \frac{1}{p}\right) = H^{1-o(1)}.$$
By repeating the proof of \Cref{higher-estimates}, one can show that
$$ \frac{1}{X} \int_X^{2X} \left(\sum_{x< n \le x+H} \mathbbm{1}_{\mathcal R(z)}(n) - R' \right)^{6} \, dx \ll H^{3+o(1)};$$
the condition $h \leq \log^{1/100} X$ ensures that the contribution of the errors from the Rosser--Iwaniec sieve are negligible\footnote{Indeed, as $P(z)$ is now much less than $X$, one could simply use the Chinese remainder theorem in place of the sieve and still obtain a negligible error term.}, and the choice $z = H^{10}$ also ensures that the error terms arising from truncating the singular series using \eqref{sdj} are also negligible for this sixth moment calculation.  We conclude that
$$ N_h(X) \ll \frac{1}{H^{6-o(1)} \log^2 X} X H^{3+o(1)}$$
giving the claim.
\end{proof}

If we sum \Cref{Nh-bound} over all $h$ and insert this into \eqref{no}, we obtain \eqref{final-bound}.  Now we turn to the conditional result \eqref{asymp}. From \eqref{no} and Lemma \ref{Nh-bound}, we see that for any fixed $H \geq 2$, one has
$$ N(X) = \sum_{2 \leq h \leq H} N_h(x) + O\left( \frac{X}{H \log^2 X} \right)$$
if $X$ is sufficiently large compared to $H$, with the implied constant here uniform in $H$.  Thus, to obtain the result, it will suffice to establish the asymptotic
\begin{equation}\label{desired}
     N_h(X) = (c_h+o(1)) \frac{X}{\log^2 X}
\end{equation}
as $X \to \infty$ for each fixed $h \geq 2$, where $c_h$ is a non-negative constant with $c_h > 0$ for at least one $h$; note from \Cref{Nh-bound} that $c_h$ is necessarily of size $O(1/h^2)$, so that the series $c \coloneqq \sum_{h=2}^\infty c_h$ is absolutely convergent to some positive quantity $c>0$.

We can restrict to even $h$, since $N_h(X)$ vanishes for odd $h$ and $X>2$; we can also exclude $h=2$ as $N_h$ also vanishes in this case.  Observe that $N_h(X)$ is precisely the set of primes $n \in [X,2X]$ such that $n, n+h$ is prime, and such that each of $n+1,\dots,n+h-1$ contain a prime factor less than $h$.  Equivalently,
$$ N_h(X) = \sum_{b \in \Omega_h} | \{ n \in [X,2X]: n = b \hbox{ mod } P(h); \quad n, n+h \hbox{ prime} \}|$$
where $\Omega_h$ is the set of residue classes $b \hbox{ mod } P(h)$ such that $b, b+h$ are coprime to $P(h)$, but $b+1,\dots,b+h-1$ are not.

We now assume the following form of the Dickson--Hardy--Littlewood prime tuples conjecture \cite{dickson,hardy-littlewood}:

\begin{conjecture}[Dickson--Hardy--Littlewood prime tuples conjecture]\label{dhltc}  Let $W \geq 1$ be fixed, and let $b_1,b_2$ be fixed distinct integers coprime to $W$.  Then, as $X \to \infty$, the number of $n \leq X$ for which $Wn+b_1$ and $Wn+b_2$ are both prime is\footnote{It would be more accurate to replace $\frac{X}{\log^2 X}$ here by $\int_2^X \frac{dt}{\log^2 t}$, but as we are content with $o()$ type error terms, we will not bother to make this refinement.}
$$ \sim \frac{X}{W \log^2 X} \prod_p \left(1 - \frac{1}{p}\right)^{-2} \left(1-\frac{\nu_p(W \cdot + b_1,W \cdot + b_2)}{p} \right),$$
where $\nu_p(W \cdot + b_1,W \cdot + b_2)$ is the number of distinct residue classes modulo $p$ found amongst $b_1, b_2$ if $p$ does not divide $W$, and zero if $p$ does divide $W$.
\end{conjecture}

\begin{table}
\begin{tabular}{|c|c|c|c|c|}
\hline
$h$ & $\Omega_h$ & $c_h$ & $\sum_{h' \leq h} c_h$ \\
\hline
4 & $\{1 \hbox{ mod } 6 \}$ & $1.3203\dots$ & $1.3203\dots$ \\
6 & $\{1,23 \hbox{ mod } 30 \}$ & $0.8802\dots$ & $2.2005\dots$ \\
8 & $\{89, 113 \hbox{ mod } 210 \}$ & $0.1760\dots$ & $2.3766\dots$ \\
10 & $\{1, 199 \hbox{ mod } 210 \}$ & $0.1760\dots$ & $2.5526\dots$ \\
12 & $\{1,199,467, 509, 1789, 1831, 2099, 2310, 2297 \hbox{ mod } 2310 \}$ & $0.0782\dots$ & $2.6308\dots$ \\
14 & 58 \hbox{ residue classes mod } 30030 & $0.0229\dots$ & $2.6824\dots$ \\
16 & 12 \hbox{ residue classes mod } 30030 & $0.0107\dots$ & $2.6931\dots$ \\
18 & 376 \hbox{ residue classes mod } 510510 & $0.0222\dots$ & $2.7154\dots$ \\
\hline
\end{tabular}
\caption{Values of $\Omega_h$ and $c_h$ for small $h$. ($\Omega_h$ is empty, and $c_h$ vanishes, when $h$ is odd or equal to two.)  The difficulty in computing $c_h$ exactly increases rapidly with $h$.}\label{ch-table}
\end{table}

Assuming this conjecture, we have for each $b \in \Omega_h$ that
\begin{align*}
 | \{ n \in [X,2X]: &n = b \hbox{ mod } P(h); \quad n, n+h \hbox{ prime} \}| \\
 &\sim \frac{X}{P(h) \log^2 X} \prod_p \left(1 - \frac{1}{p}\right)^{-2} \left(1-\frac{\nu_p(W \cdot + b,W \cdot + b+h)}{p} \right) \\
 &\sim \frac{X {\mathfrak S}}{\left( \prod_{2 < p < h} (p-2) \right) \log^2 X},
\end{align*}
where the twin prime constant ${\mathfrak S}$ was defined in \eqref{twinprime}, and where we have used the fact that $h$ is even but not equal to $2$.  We thus obtain the desired asymptotic \eqref{desired} with the explicit constant
$$ c_h \coloneqq {\mathfrak S} \frac{|\Omega_h|}{\prod_{2 < p < h} (p-2)}.$$
As $c = \sum_{h=2}^\infty c_h$, to conclude the theorem it suffices to show that $c_h$ is positive for at least one $h$.  But the example provided by Erd\H{o}s already suffices for this; explicitly, we can take $h=18$, so that $P(18) = 510510$, and then $\Omega_{18}$ contains the residue class $2183 \hbox{ mod } 510510$ (in fact, as already noted by Erd\H{o}s, it contains the larger residue class $2183 \hbox{ mod } 30030$).  There are simpler examples; for instance, $\Omega_4$ contains the residue class $1 \hbox{ mod } P(4) = 1 \hbox{ mod } 6$, reflecting the fact that any prime gap of the form $(6n+1, 6n+5)$ will not contain a rough number.
We give the values of $\Omega_h$ and $c_h$ for small $h$ in \Cref{ch-table}; as we expect $c_h$ to decay rapidly, we therefore expect $c \in [2.7, 2.8]$, although we did not rigorously establish the upper bound.

\begin{remark}  As the bound $c_h \ll 1/h^2$ already indicates, the set $\Omega_h$ of residue classes should become extremely sparse as $h$ increases; but it should be possible to adapt the constructions of large prime gaps of Rankin \cite{rankin} and subsequent authors (e.g., \cite{fgkmt}) to show that $\Omega_h$ remains non-empty (or equivalently, $c_h$ non-zero) for sufficiently large even $h$. We will not pursue this question further here.
\end{remark}

\bibliographystyle{amsplain}

\end{document}